\documentclass[11pt]{article}
\usepackage{amsmath, amsthm, amscd, amsfonts, amssymb}
\date{ }

\newcommand{\al}{\alpha}

\newtheorem{theorem}{Theorem}[section]

\newtheorem{corollary}[theorem]{Corollary}
\newtheorem{definition}[theorem]{Definition}
\newtheorem{remark}[theorem]{Remark}
\newtheorem{conjecture}[theorem]{Conjecture}
\newtheorem{problem}[theorem]{Problem}
\newtheorem{main}[theorem]{Main Theorem}

\title{\bf The Answers to a Problem and Two Conjectures about OD-Characterization of Finite Groups }
\smallskip
\normalsize
\author{{\bf Ali Mahmoudifar  \&  Behrooz Khosravi}
\\ Dept. of Pure  Math.,  Faculty  of Math. and Computer Sci.,\\
Amirkabir University of Technology (Tehran Polytechnic),\\
424, Hafez Ave., Tehran 15914, IRAN, \\  e-mail:
khosravibbb@yahoo.com\\alimahmoudifar@gmail.com}
\begin{document}
\maketitle
\begin{abstract}
In [Akbari and Moghaddamfar, Recognizing by order and degree pattern of some projective special linear groups, {\it Internat. J. Algebra Comput.}, 2012]
the authors possed the following problem:
\\
{\bf Problem.} {\it Is there a simple group which is $k$-fold OD-characterizable for $k\geq3\ ?$ }

In this paper as the main result we give positive answer to the above problem and we introduce
two simple groups which are $k$-fold OD-characterizable such that $k\geq6$.

Also in [R. Kogani-Moghadam and A. R. Moghaddamfar, Groups with the same order and degree pattern,
{\it Science China Mathematics}, 2012], the authors possed two conjectures as follows:
\\
{\bf Conjecture 1.} {\it All alternating groups $A_m$ with $m \not= 10$ are OD-characterizable.}
\\
{\bf Conjecture 2.} {\it All symmetric groups $S_m$, with $m \not= 10$, are $n$-fold OD-characterizable, where
$n\in\{1, 3\}$.}

In this paper we find some alternating and some symmetric groups such that these conjectures are not true for them.
\end{abstract}
\section{Introduction}

In this paper every group is finite. If $n$ is a natural number, then we
denote by $\pi(n)$, the set of all prime divisors of $n$. If $G$ is a finite group, then $\pi(|G|)$ is denoted by $\pi(G)$.
The prime graph $GK(G)$ of a group $G$ is defined as a graph with vertex set $\pi(G)$ in
which two distinct primes $p,q\in\pi(G)$ are adjacent (and we write $p\sim q$) if G contains an element of order
$pq$. In {\rm\cite[Proposition~1.1]{anadja}}\label{adalt}, it is proved that if $r, s\in\pi(A_n)\setminus\{2\}$, then $r\nsim s$ in $GK(A_n)$ if and only if $r + s > n$ and
$2\nsim r$ in $GK(A_n)$ if and only if $r + 4 > n$.
\begin{definition}{\rm(\cite{mogh})}
Let $G$ be a group and $\pi(G)=\{p_1,p_2,\ldots,p_k\}$ where $p_1<p_2<\cdots<p_k$.
Then the degree pattern of $G$ is defined as follows:
$$D(G):=(\deg(p_1),\deg(p_2),\ldots,\deg(p_k)),$$
where $\deg(p_i)$, $1\leq i\leq k$, is the degree of vertex $p_i$ in $GK(G)$.
\end{definition}
For a group $G$, denote by $h_{OD}(G)$ the number of nonisomorphism finite groups $H$, such that $|H| = |G|$
and $D(H) = D(G)$.
\begin{definition}{\rm(\cite{mogh})}
Let $G$ be a group. We call that $G$ is $k$-fold OD-characterizable if $h_{OD}(G)=k$.
When $k=1$, the finite group $G$ is called OD-characterizable.
\end{definition}
OD-characterization of some finite groups are considered by some authors (see the references of \cite{mogh}). In \cite{hads}, the following problem
about OD-characterization of finite simple groups is possed (see also \cite{in,mogh,hads}):
 \begin{problem}\label{problem}
 Is there a simple group which is $k$-fold OD-characterizable for $k\geq3\ ?$
\end{problem}

Also in \cite{mogh}, the authors put the following conjectures:
\begin{conjecture}\label{conjecture1}
All alternating groups $A_m$ with $m\not= 10$ are OD-characterizable.
\end{conjecture}
\begin{conjecture}\label{conjecture2}
All symmetric groups $S_m$, with $m\not= 10$, are $n$-fold OD-characterizable, where
$n \in\{1,3\}$.
\end{conjecture}

Until now, Conjecture 1.4 is valid for $n=p$, $p+1$ and $p+2$,  where $p$ is an odd prime number.

In this paper as the main result we prove the following theorem and using this theorem we introduce some finite simple groups which are $k$-fold OD-characterizable, where $k\geq 6$. Also
we present some counterexamples for the above conjectures.
\begin{main}\label{main}
Let $p$ be an odd prime number such that $p+2$ and $p+4$ are not prime and $p+6=5^{\al}$, for some natural number $\al$. Also let
$H$ and $T$ be some groups such that $|H|=p+6=5^{\al}$ and $|T|=2(p+6)$. Then we have:
\begin{enumerate}
\item[{\rm (1)}] $|A_{p+6}|=|A_{p+5}\times H|$ and $D(A_{p+6})=D(A_{p+5}\times H)$,

\item[{\rm (2)}] $|S_{p+6}|=|S_{p+5}\times H|$ and $D(S_{p+6})=D(S_{p+5}\times H)$,

\item[{\rm (3)}]   $|S_{p+6}|=|A_{p+5}\times T|$ and $D(S_{p+6})=D(A_{p+5}\times T)$.
\end{enumerate}
\end{main}

\begin{remark}
Easily we can see that $p=619=5^4-6$, $15619=5^6-6$, $9765619=5^{10}-6$ satisfy the hypothesis of the main theorem.
\end{remark}

\begin{corollary}
If $(G,K)=(A_{625},A_{624})$ or $(S_{625},S_{624})$, then the following groups have the same order and degree pattern:
$G$, $K\times Z_{625}$, $K\times Z_{125}\times Z_{5}$, $K\times Z_{25}\times Z_{25}$, $K\times Z_{25}\times Z_5\times Z_5$, $K\times Z_5 \times Z_5 \times Z_5 \times Z_5$.
\end{corollary}
\begin{corollary}
We have $h_{OD}(A_{625})\geq6$. In Particular, the simple group $A_{625}$ is an answer for Problem{\rm~\ref{problem}} and also we see that Conjecture{\rm~\ref{conjecture1}} is not true.
\end{corollary}
\begin{corollary}
We have $h_{OD}(S_{625})\geq6$. In Particular, Conjecture{\rm~\ref{conjecture2}} is not true.
\end{corollary}
\begin{corollary}
Let $p=5^{\al}-6$ satisfies the assumptions of the main theorem. If $m$ is the number of nonisomorphic finite groups of order
$5^{\al}$, then $h_{OD}(A_{p+6})\geq m+1$ and $h_{OD}(S_{p+6})\geq m+1$.
\end{corollary}
\section{Proof of The Main Theorem}
Let $p$ be a prime number such that $p+2$ and $p+4$ are not prime and $p+6=5^{\al}$, for some natural number $\al$.
The proof of the main theorem for symmetric groups is similar to alternating groups.

First we show that $GK(A_{p+6})=GK(A_{p+5})$. By assumption,  $p+6$ is a power of $5$ and so $\pi(A_{p+6})=\pi(A_{p+5})$. Since $p+2$ and $p+4$ are not prime numbers, by \cite[Proposition~1.1]{anadja}, we get that  every vertex is adjacent to $2$, $3$ and $5$ in $GK(A_{p+5})$ and $GK(A_{p+6})$. Then it is sufficient to check the adjacency of odd primes in these graphs.

Let $r,s\in \pi(A_{p+5})\setminus\{2\}$. By \cite[Proposition~1.1]{anadja}, if $r\sim s$ in $GK(A_{p+5})$, then $r+s\leq p+5<p+6$ and so $r\sim s$ in $GK(A_{p+6})$. Let $r\nsim s$, in $GK(A_{p+5})$. Then by \cite[Proposition~1.1]{anadja}, we have $r+s>p+5$. Since $p+6$ is odd, we get that $r+s\not=p+6$, and so $r+s>p+6$. This implies that $r\nsim s$ in $GK(A_{p+6})$.

Let $r,s\in\pi(A_{p+6})\setminus\{2\}$. If $r\nsim s$ in $GK(A_{p+6})$, then by \cite[Proposition~1.1]{anadja}, $r+s>p+6>p+5$ and so $r\nsim s$ in $GK(A_{p+5})$. Since $r+s\not=p+6$, if $r\sim s$ in $GK(A_{p+6})$, then $r+s\leq p+5$ and so $r\sim s$ in $GK(A_{p+5})$.

Therefore by the above discussion we conclude that $GK(A_{p+6})=GK(A_{p+5})$. Now if $H$ is a group such that $|H|=p+6=5^{\al}$, then $|A_{p+6}|=|A_{p+5}\times H|$. Also since $5$ is adjacent to each vertex in $GK(A_{p+6})$ and $|H|=5^{\al}$, then we get that $GK(A_{p+6})=GK(A_{p+5}\times H)$.
Therefore for every group $H$, such that $|H|=p+6=5^{\al}$, we have $|A_{p+6}|=|A_{p+5}\times H|$ and $D(A_{p+6})=D(A_{p+5}\times H)$.


\end{document}